\documentclass[final,3p,times,twocolumn]{elsarticle}
\usepackage{amsfonts}

\usepackage{eurosym}
\usepackage{amssymb}
\usepackage{amsthm}
\usepackage{amsmath}
\usepackage{adjustbox}

\journal{Physica D}

\begin{document}

\begin{frontmatter}

\title{Differentiable programming and its applications to dynamical systems}

\author{Adri\'{a}n Hern\'{a}ndez and Jos\'{e} M. Amig\'{o}\corref{cor1}}

\address{Centro de Investigaci\'{o}n Operativa, Universidad Miguel Hern\'{a}ndez, Avenida de la Universidad s/n, 03202 Elche, Spain}
\cortext[cor1]{Corresponding author}

\begin{abstract}
Differentiable programming is the combination of classical neural networks modules with algorithmic ones in an end-to-end differentiable model. These new models, that use automatic differentiation to calculate gradients, have new learning capabilities (reasoning, attention and memory). In this tutorial, aimed at researchers in nonlinear systems with prior knowledge of deep learning, we present this new programming paradigm, describe some of its new features such as attention mechanisms, and highlight the benefits they bring. Then, we analyse the uses and limitations of traditional deep learning models in the modeling and prediction of dynamical systems. Here, a dynamical system is meant to be a set of state variables that evolve in time under general internal and external interactions. Finally, we review the advantages and applications of differentiable programming to dynamical systems.
\end{abstract}

\begin{keyword}
Deep learning \sep differentiable programming \sep dynamical systems \sep attention \sep recurrent neural networks
\end{keyword}

\end{frontmatter}



\section{Introduction}

\label{sec1}

The increase in computing capabilities together with new deep learning
models has led to great advances in several machine learning tasks \cite%
{LeCun2015DeepLearning, Sutskever2014SequenceTS, Silver2017MasteringTG}.

Deep learning architectures such as Recurrent Neural Networks (RNNs) and
Convolutional Neural Networks (CNNs), as well as the use of distributed
representations in natural language processing, have allowed to take into
account the symmetries and the structure of the problem to be solved.

However, a major criticism of deep learning remains, namely, that it only
performs perception, mapping inputs to outputs \cite{Marcus2018Deep}.

A new direction to more general and flexible models is differentiable
programming, that is, the combination of geometric modules (traditional
neural networks) with more algorithmic modules in an end-to-end
differentiable model. As a result, differentiable programming is a dynamic computational graph composed of differentiable functions that provides
not only perception but also reasoning, attention and memory. To efficiently calculate derivatives, this approach uses automatic differentiation, an algorithmic technique similar to backpropagation and implemented in modern software packages such as PyTorch, Julia, etc.

To keep our exposition concise, this tutorial is aimed at researchers in nonlinear systems with prior knowledge of deep learning; see \cite{Goodfellow2016DL} for an excellent introduction to the concepts and methods of deep learning. Therefore, this tutorial focuses right away on the limitations of traditional deep learning and the advantages of differential programming, with special attention to its application to dynamical systems. By a dynamical system we mean here and hereafter a set of state variables that evolve in time under the influence of internal and possibly external inputs.

Examples of differentiable programming techniques that have been successfully developed in recent years include 

(i) attention mechanisms \cite{Bahdanau2014Attention}, which allow the model
to automatically search and learn which parts of a source sequence are relevant to predict the target element,

(ii) self-attention, 

(iii) end-to-end Memory Networks \cite{Suk2015EndToEndMN}, and 

(iv) Differentiable Neural Computers (DNCs) \cite{Graves2016DNC}, which are
neural networks (controllers) with an external read-write memory.

As expected, in recent years there has been a growing interest in applying
deep learning techniques to dynamical systems. In this regard, RNNs and Long
Short-Term Memories (LSTMs), specially designed for sequence modelling and temporal dependence, have been successful in various applications to dynamical
systems such as model identification and time series prediction \cite{WangModelIdLSTM, Yu2017ConceptLSTM, Li2018LSTMTourismFlow}.

The performance of theses models (e.g. encoder-decoder networks), however,
degrades rapidly as the length of the input sequence increases and they are
not able to capture the dynamic (i.e., time-changing) interdependence between time steps. The combination of neural networks with new differentiable modules could overcome some of those problems and offer new opportunities and applications.

Among the potential applications of differentiable programming to dynamical
systems let us mention 

(i) attention mechanisms to select the relevant time
steps and inputs, 

(ii) memory networks to store historical data from dynamical
systems and selectively use it for modelling and prediction, and 

(iii) the use of differentiable components in scientific computing. 

\noindent Despite some achievements, more work is still needed to verify the benefits of these models over traditional networks.

Thanks to software libraries that facilitate automatic differentiation,
differentiable programming extends deep learning models with new
capabilities (reasoning, memory, attention, etc.) and the models can be efficiently coded and implemented.

In the following sections of this tutorial we introduce differentiable programming and explain in detail why it is an extension of deep learning (Section 2). We describe some models based on this new approach such as attention mechanisms (Section 3.1), memory networks and differentiable neural computers (Section 3.2), and continuous learning (Section 3.3). Then we review the use of deep learning in dynamical systems and their limitations (Section 4.1). And, finally, we present the new opportunities that differentiable programming can bring to the modelling, simulation and prediction of dynamical systems (Section 4.2). The conclusions and outlook are summarized in Section 5.

\section{From deep learning to differentiable programming}

\label{sec2}

In recent years, we have seen major advances in the field of machine
learning. The combination of deep neural networks with the computational
capabilities of Graphics Processing Units (GPUs) \cite{Yadan2013MultiGPU}
has improved the performance of several tasks (image recognition, machine
translation, language modelling, time series prediction, game playing and more) 
\cite{LeCun2015DeepLearning, Sutskever2014SequenceTS, Silver2017MasteringTG}. Interestingly, deep learning models and architectures have evolved to take into account the structure of the problem to be resolved.

\textit{Deep learning} is a part of machine learning that is based on neural networks and uses multiple layers, where each layer extracts higher level features from the input. RNNs are a special class of neural networks where outputs from previous steps are fed as inputs to the current step \cite%
{Graves2009ANC,Sherstinsky2018FundRNN}. This recurrence makes them
appropriate for modelling dynamic processes and systems.

CNNs are neural networks that alternate convolutional and pooling layers to implement translational invariance \cite%
{Lecun1998CNN}. They learn spatial hierarchies of features through
backpropagation by using these building layers. CNNs are being applied successfully
to computer vision and image processing \cite{Yamashita2018CNN}.

Especially important is the use of distributed representations as inputs to
natural language processing pipelines. With this technique, the words of the
vocabulary are mapped to an element of a vector space with a much lower
dimensionality \cite{Bengio2003NeuralModel, Mikolov2013Word2vec}. This word
embedding is able to keep, in the learned vector space, some of the
syntactic and semantic relationships presented in the original data.

\begin{figure}[h]
\centering
\includegraphics[scale=0.58]{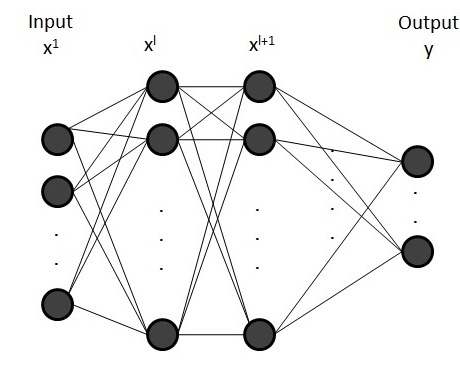} 
\caption{Multilayer neural network.}
\label{mnn}
\end{figure}

Let us recall that, in a feedforward neural network (FNN) composed of multiple layers, the output (without the bias term) at layer $l$, see Figure \ref{mnn}, is defined as

\begin{equation}  \label{eq:mnn}
\boldsymbol{x}^{l+1}=\sigma(W^{l}\boldsymbol{x}^{l}),
\end{equation}

\noindent $W^{l}$ being the weight matrix at layer $l$. $\sigma$ is the activation function and $\boldsymbol{x}^{l+1}$, the output vector at layer $l$ and the input vector at layer $l+1$. The weight matrices for the different layers are the parameters of the model.

Learning is the mechanism by which the parameters of a neural network are
adapted to the environment in the training process. This is an optimization
problem which has been addressed by using gradient-based methods, in which
given a cost function $f:\mathbb{R}^n\rightarrow \mathbb{R}$, the algorithm finds local minima $w^{*}=\operatorname{arg\,min}_{w} f(w)$ updating each layer parameter $w_{ij}$ with the rule $w_{ij}:=w_{ij}-\eta \nabla_{w_{ij}} f(w)$, $\eta>0$ being the learning rate.

In addition to regarding neural networks as universal approximators, there is no sound theoretical explanation for a good performance of deep learning. Several theoretical frameworks have been proposed:
\begin{enumerate}[(i)]

\item As pointed out in \cite{Lin2016Cheaplearning}, the class of functions of practical interest can be approximated with exponentially fewer parameters than the generic ones. Symmetry, locality and compositionality properties make it possible to have simpler neural networks.

\item From the point of view of information theory \cite{ShwartzZiv2017ITDL}, an explanation has been put forward based on how much information each layer of the neural network retains and how this information varies with the training and testing process.

\end{enumerate} 

Although deep learning can implicitly implement logical reasoning \cite%
{Hohenecker2018Onto}, it has limitations that make it difficult to achieve more general intelligence \cite{Marcus2018Deep}. Among these limitations, we can highlight the following:
\begin{enumerate}[(i)]
\item It only performs perception, representing a mapping between inputs and outputs.
\item It follows a hybrid model where synaptic weights perform both processing and memory tasks but doesn't have an explicit external memory.
\item It does not carry out conscious and sequential reasoning, a process that is based on perception and memory through attention.

\end{enumerate} 

A path to a more general intelligence, as we will see below,
is the combination of geometric modules with more algorithmic modules in an
end-to-end differentiable model. This approach, called differentiable
programming, adds new parametrizable and differentiable components to
traditional neural networks.

\textit{Differentiable programming}, a broad term, is defined in \cite{Wang2018Backprop} as a programming model (model of how a computer program is executed), trainable with gradient descent, where neural networks are truly functional blocks with data-dependent branches and recursion.

Here, and for the purposes of this tutorial, we define differentiable programming as a programming model with the following characteristics:
\begin{enumerate}[(i)]
\item Programs are directed acyclic graphs.

\item Graph nodes are mathematical functions or variables and the edges correspond to the flow of intermediate values between the nodes.

\item $n$ is the number of nodes and $l$ the number of input variables of the graph, with $1\leq l< n$. $v_i$ for $i \in \{1,...,n\}$ is the variable associated with node $i$.

\item $E$ is the set of edges in the graph. For each $(i,j) \in E$ we have $i<j$, therefore the graph is topologically ordered.

\item $f_i$ for $i \in \{(l+1),...,n\}$ is the differentiable function computed by node $i$ in the graph. $\alpha_i$ for $i \in \{(l+1),...,n\}$ contains all input values for node $i$.

\item The forward algorithm or pass, given input variables $v_1,...,v_l$ calculates $v_i=f_i(\alpha_i)$ for $i=\{(l+1),...,n\}$.

\item The graph is dynamically constructed and composed of parametrizable functions that are differentiable and whose parameters are learned from data.

\end{enumerate}

Then, neural networks are just a class of these differentiable programs composed of classical blocks (feedforward, recurrent neural networks, etc.) and new ones such as differentiable branching, attention, memories, etc.

Differentiable programming can be seen as a continuation of the deep learning
end-to-end architectures that have replaced, for example, the traditional
linguistic components in natural language processing \cite{Deng2018NLP,
Goldberg2017NLPBook}. To efficiently calculate the derivatives in a gradient
descent, this approach uses automatic differentiation, an algorithmic
technique similar but more general than backpropagation.

Automatic differentiation, in its reverse mode and in contrast to manual,
symbolic and numerical differentiation, computes the derivatives in a
two-step process \cite{Baydin2018AutomaticDiff, Wang2018DemystifyingDP}. 
As described in \cite{Baydin2018AutomaticDiff} and rearranging the indexes of the previous definition, a function $f:R^{n}\rightarrow R^{m}$ is constructed with intermediate variables $v_{i}$ such that:
\begin{enumerate}[(i)]
\item variables $v_{i-n}=x_{i}, i=1,...,n$ are the inputs variables.
\item variables $v_{i}, i=1,...,l$ are the intermediate variables.
\item variables $y_{m-i}=v_{l-i}, i=m-1,...,0$ are the output variables.
\end{enumerate}

In a first step, similar to the forward pass described before, the computational graph is built populating intermediate variables $v_i$ and recording the dependencies. In a second step, called the backward pass, derivatives are calculated by propagating for the output $y_j$ being considered, the adjoints $\overline{v}_{i}=\frac{\partial {y}_j}{\partial {v}_{i}}$ from the output to the inputs.

The reverse mode is more efficient to evaluate for functions with a large number of inputs (parameters) and a small number of outputs. When $f:R^{n}\rightarrow R$, as is the case in machine learning with $n$ very large and $f$ the cost function, only one pass of the reverse mode is necessary to compute the gradient $\nabla f=(\frac{\partial y}{\partial x_1},...,\frac{\partial y}{\partial x_n}).$ 
 
In the last years, deep learning frameworks such as PyTorch
have been developed that provide reverse-mode automatic differentiation \cite%
{Paszke2017automatic}. The define-by-run philosophy of PyTorch, whose
execution dynamically constructs the computational graph, facilitates the development of general differentiable programs.

Differentiable programming is an evolution of classical (traditional) software programming where, as shown in Table 1:
\begin{enumerate}[(i)]
\item Instead of specifying explicit instructions to the computer, an objective is set and an optimizable architecture is defined which allows to search in a subset of possible programs.
\item The program is defined by the input-output data and not predefined by the user.
\item The algorithmic elements of the program have to be differentiable, say, by converting them into differentiable blocks. 
\end{enumerate}

\begin{table}[h!]
\begin{center}
\begin{adjustbox}{width=\columnwidth,center}
    \begin{tabular}{l|l|l} 
      \textbf{Classical} & \textbf{Differentiable}\\
      \hline
      Sequence of instructions & Sequence of diff. primitives\\
      Fixed architecture & Optimizable architecture\\
      User defined & Data defined\\
      Imperative programming & Declarative programming\\
      Intuitive & Abstract\\   
    \end{tabular}
    \end{adjustbox} 
\end{center}
\caption{Differentiable vs classical programming.}
\label{table1}
\end{table}

RNNs, for example, are an evolution of feedforward networks because they are
classical neural networks inside a for-loop (a control flow statement for
iteration) which allows the neural network to be executed repeatedly with
recurrence. However, this for-loop is a predefined feature of the model. Differentiable programming allows to dynamically constructs the graph and vary the length of the loop. Then, the ideal situation would be to augment the neural network with programming primitives (for-loops, if branches, while statements, external memories, logical modules, etc.) that are not predefined by the user but are parametrizable by the training data.

\begin{figure}[h]
\centering
\includegraphics[scale=0.60]{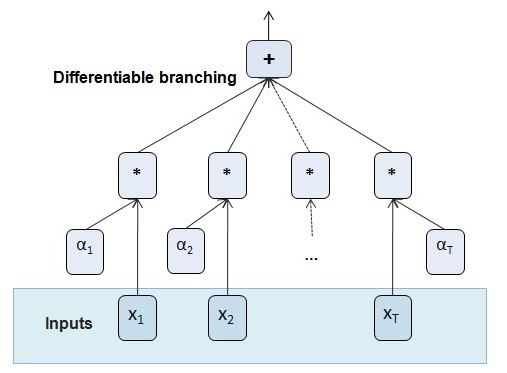} 
\caption{Computational graph of differentiable branching.}
\label{DiffBranching}
\end{figure}

The trouble is that many of these programming primitives are not
differentiable and need to be converted into optimizable modules. For
instance, if the condition $a$ of an "if" primitive (e.g., if $a$ is satisfied do $y(x)$, otherwise do $z(x)$) is to be learned, it can be the output of a
neural network (linear transformation and a sigmoid function) and the
conditional primitive will transform into a weighted combination of both
branches $ay(x)+(1-a)z(x)$. Similarly, in an attention module, different
weights that are learned with the model are assigned to give a different
influence to each part of the input. Figure \ref{DiffBranching} shows the computational graph of a conditional branching.

The process of extending deep learning with differentiable primitives would consist of the following steps:
\begin{enumerate}[(i)]

\item Select a new function that improves the classical input-output transformation of deep learning, e.g. attention, continuous learning, memories, etc.
\item Convert this function into a directed acyclic graph, a sequence of parametrizable and differentiable functions. For example, Figure \ref{DiffBranching} shows this sequence of operations used in attention for differentiable branching. 
\item Integrate this new function into the base model.

\end{enumerate}

In this way, using differentiable programming we can combine traditional
perception modules (CNN, RNN, FNN) with additional algorithmic modules that
provide reasoning, abstraction and memory \cite{Yang2017DiffeLogical}. In
the following section we describe, by following this process, some examples of this approach that have been developed in recent years.

\section{Differentiable learning and reasoning}

\label{sec3}

\subsection{Differentiable attention}

\label{sec31}

One of the aforementioned limitations of deep learning models is that they do not perform conscious and sequential reasoning, a process that is based on perception and memory through attention. 

Reasoning is the process of consciously establishing and verifying facts combining attention with new or existing information. An attention mechanism allows the brain to focus on one part of the input or memory (image, text, etc), giving less attention to others.

Attention mechanisms have provided and will provide a paradigm shift in machine learning. From traditional large-scale vector transformations to a more conscious process that focuses only on a set of elements, e.g. decomposing a problem into a sequence of attention based reasoning operations \cite{Hudson2018CompAttention}.

\begin{figure}[h]
\centering
\includegraphics[scale=0.50]{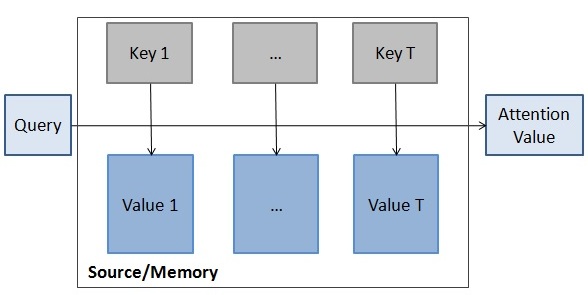} 
\caption{Attention diagram.}
\label{KeyValue}
\end{figure}

One way to make this attention process differentiable is to make it a convex combination of the input or memory, where all the steps are differentiable and the combination weights are parametrizable.

As in \cite{Vaswani2017AttentionNeed}, this differentiable attention process is described as mapping a query and a set of key-value pairs to an output:

\begin{equation}  \label{eq:attint1}
att(\boldsymbol{q},\boldsymbol{s})=\sum_{i=1}^{T}\alpha_{i}(\boldsymbol{q},\boldsymbol{k}_{i})\boldsymbol{V}_{i},
\end{equation} 

\noindent where, as seen in figure \ref{KeyValue}, $\boldsymbol{k}_{i}$ and $\boldsymbol{V}_{i}$ are the key and the value vectors from the source/memory $\boldsymbol{s}$, and $\boldsymbol{q}$ is the query vector. $\alpha_{i}(\boldsymbol{q},\boldsymbol{k}_{i})$ is the similarity function between the query and the corresponding key and is calculated by applying the softmax function: 

\begin{equation}  \label{eq:softmax}
Softmax(z_i)=\frac{exp(z_i)}{\sum_{i^{\prime}}exp(z_{i^\prime})} 
\end{equation}

\noindent to the score function $score(\boldsymbol{q},\boldsymbol{k}_{i}):$

\begin{equation}  \label{eq:attint2}
\alpha_{i}=\frac{\exp(score(\boldsymbol{q},\boldsymbol{k}_{i}))}{\sum_{i^{\prime}=1}^{T}\exp(score(\boldsymbol{q},\boldsymbol{k}_{i^{\prime}}))}.
\end{equation}

The score function can be computed using a feedforward neural network:
 
\begin{equation}  \label{eq:attint3}
score(\boldsymbol{q},\boldsymbol{k}_{i})=\boldsymbol{Z}_{a}\tanh(\boldsymbol{W}_{a}[\boldsymbol{q},\boldsymbol{k}_{i}])),
\end{equation}

\noindent as proposed in \cite{Bahdanau2014Attention}, where $\boldsymbol{Z}_{a}$ and $\boldsymbol{W}_{a}$ are matrices to be jointly learned with the rest of
the model and $[\boldsymbol{q},\boldsymbol{k}_{i}]$ is a linear function or concatenation of $\boldsymbol{q}$ and $\boldsymbol{k}_{i}$. Also, in \cite{Graves2014NTM} the authors use a cosine similarity measure for content-based attention, namely,

\begin{equation}  \label{eq:attint4}
score(\boldsymbol{q},\boldsymbol{k}_{i})=\cos((\boldsymbol{q},\boldsymbol{k}_{i}))
\end{equation}

\noindent where $((\boldsymbol{q},\boldsymbol{k}_{i}))$ denotes the angle between $\boldsymbol{q}$ and $\boldsymbol{k}_{i}$.

Then, differentiable attention can be seen as a sequential process of reasoning in which the task (query) is guided by a set of elements of the input source (or memory) using attention. 

The attention process can focus on:
\begin{enumerate}[(i)]
\item Temporal dimensions, e.g. different time steps of a sequence.
\item Spatial dimensions, e.g. different regions of an image.
\item Different elements of a memory.
\item Different features or dimensions of an input vector, etc. 

\end{enumerate}

Depending on where the process is initiated, we have:
\begin{enumerate}[(i)]

\item Top-down attention, initiated by the current task.
\item Bottom-up, initiated spontaneously by the source or memory.

\end{enumerate}

\subsubsection{Attention mechanisms in seq2seq models}
\label{sec311}
\begin{figure}[h]
\centering
\includegraphics[scale=0.6]{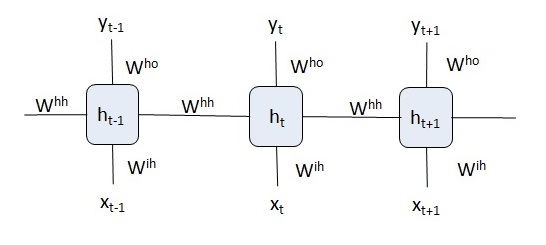} 
\caption{Temporal structure of a recurrent neural network.}
\label{rnn}
\end{figure}

RNNs (see Figure \ref{rnn}) are a basic component of modern deep learning architectures, especially of encoder-decoder networks. The following equations define the time evolution of an RNN: 
\begin{equation}  \label{eq:rnn1}
\boldsymbol{h}_{t}=f^{h}(W^{ih}\boldsymbol{x}_{t}+W^{hh}\boldsymbol{h}%
_{t-1}),
\end{equation}
\begin{equation}  \label{eq:rnn2}
\boldsymbol{y}_{t}=f^{o}(W^{ho}\boldsymbol{h}_{t}),
\end{equation}

\noindent $W^{ih}$, $W^{hh}$ and $W^{ho}$ being weight matrices. $f^{h}$ and $f^{o}$
are the hidden and output activation functions while $\boldsymbol{x}_{t}$, $%
\boldsymbol{h}_{t}$ and $\boldsymbol{y}_{t}$ are the network input, hidden state
and output.

An evolution of RNNs are LSTMs \cite{Hochreiter1997LSTM}, an RNN structure with gated units, i.e. regulators. LSTM are composed of a cell, an input gate, an output gate and a forget gate, and allow gradients to flow unchanged. The memory cell remembers values over arbitrary time intervals and the three gates regulate the flow of information into and out of the cell.

An encoder-decoder network maps an input sequence to a target one with both sequences of arbitrary length \cite{Sutskever2014SequenceTS}. They have
applications ranging from machine translation to time series prediction.

\begin{figure}[h]
\centering
\includegraphics[scale=0.6]{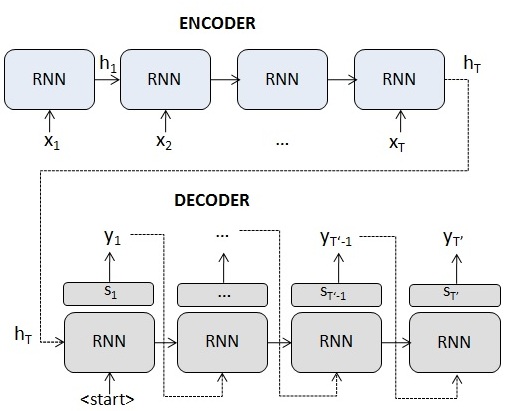} 
\caption{An encoder-decoder network.}
\label{encoder}
\end{figure}

More specifically, this mechanism uses an RNN (or any of its variants, an LSTM or a GRU, Gated Recurrent Unit) to map the input sequence to
a fixed-length vector, and another RNN (or any of its variants) to decode the target sequence
from that vector (see Figure \ref{encoder}). Such a seq2seq model features normally an architecture composed of:

\begin{enumerate}[(i)]

\item An encoder which, given an input sequence $X=(\boldsymbol{x}_{1},%
\boldsymbol{x}_{2},...,\boldsymbol{x}_{T})$ with $\boldsymbol{x}_{t}\in 
\mathbb{R}^n$, maps $\boldsymbol{x}_{t}$ to

\begin{equation}  \label{eq:enc}
\boldsymbol{h}_{t}=f_{1}(\boldsymbol{h}_{t-1},\boldsymbol{x}_{t}),
\end{equation}

where $\boldsymbol{h}_{t} \in \mathbb{R}^m$ is the hidden state of the
encoder at time $t$, $m$ is the size of the hidden state and $f_{1}$ is an
RNN (or any of its variants).

\item A decoder, where $\boldsymbol{s}_{t}$ is the hidden state and whose initial state $\boldsymbol{s}_{0}$ is initialized
with the last hidden state of the encoder $\boldsymbol{h}_{T}$. It generates
the output sequence $Y=(\boldsymbol{y}_{1},\boldsymbol{y}_{2},...,%
\boldsymbol{y}_{T^{\prime }})$, $\boldsymbol{y}_{t}\in \mathbb{R}^o$ (the dimension $o$ depending on the task), with

\begin{equation}  \label{eq:dec}
\boldsymbol{y}_{t}=f_{2}(\boldsymbol{s}_{t-1},\boldsymbol{y}_{t-1}),
\end{equation}

where $f_{2}$ is an RNN (or any of its variants) with an additional softmax layer.

\end{enumerate}

Because the encoder compresses all the information of the input sequence in
a fixed-length vector (the final hidden state $\boldsymbol{h}_{T}$), the
decoder possibly does not take into account the first elements of the input
sequence. The use of this fixed-length vector is a limitation to improve the
performance of the encoder-decoder networks. Moreover, the performance 
of encoder-decoder networks degrades rapidly as the
length of the input sequence increases \cite{cho2014ProblemsEncDecod}. This occurs  in applications such as machine translation and time series predition, where  it is necessary to model long time dependencies.

The key to solve this problem is to use an attention mechanism. In \cite{Bahdanau2014Attention} an extension of the basic encoder-decoder arquitecture was proposed by allowing the model to automatically search and
learn which parts of a source sequence are relevant to predict the target
element. Instead of encoding the input sequence in a fixed-length vector, it
generates a sequence of vectors, choosing the most appropriate subset of
these vectors during the decoding process.

With the attention mechanism, the encoder is a bidirectional RNN \cite%
{Graves2013HybridRNN} with a forward hidden state $\overrightarrow{%
\boldsymbol{h}_{i}}=f_{1}(\overrightarrow{\boldsymbol{h}}_{i-1},\boldsymbol{x%
}_{i})$ and a backward one $\overleftarrow{\boldsymbol{h}_{i}}=f_{1}(%
\overleftarrow{\boldsymbol{h}}_{i+1},\boldsymbol{x}_{i})$. The encoder state
is represented as a simple concatenation of the two states,

\begin{equation}  \label{eq:bidirectional}
\boldsymbol{h}_{i}=[{\overrightarrow{\boldsymbol{h}_{i}}};{\overleftarrow{%
\boldsymbol{h}_{i}}}],
\end{equation}

\noindent with $i=1,...,T$. The encoder state includes
both the preceding and following elements of the sequence, thus capturing
information from neighbouring inputs.

The decoder has an output

\begin{equation}  \label{eq:decatt}
\boldsymbol{y}_{t}=f_{2}(\boldsymbol{s}_{t-1},\boldsymbol{y}_{t-1},%
\boldsymbol{c}_{t})
\end{equation}

\noindent for $t=1,...,T^{\prime }$. $f_{2}$ is an RNN with an additional softmax layer, and the input is a concatenation of $\boldsymbol{y}_{t-1}$ with the context vector $\boldsymbol{c}_{t}$, which is a sum of hidden states of the
input sequence weighted by alignment scores:

\begin{equation}  \label{eq:att1}
\boldsymbol{c}_{t}=\sum_{i=1}^{T}\alpha_{ti}\boldsymbol{h}_{i}.
\end{equation}

\noindent Similar to equation (\ref{eq:attint2}), the weight $\alpha_{ti}$ of each state $\boldsymbol{h}_{i}$ is calculated by

\begin{equation}  \label{eq:att2}
\alpha_{ti}=\frac{\exp(score(\boldsymbol{s}_{t-1},\boldsymbol{h}_{i}))}{%
\sum_{i^{\prime }=1}^{T}\exp(score(\boldsymbol{s}_{t-1},\boldsymbol{h}%
_{i^{\prime}}))}.
\end{equation}

\begin{figure}[h]
\centering
\includegraphics[scale=0.6]{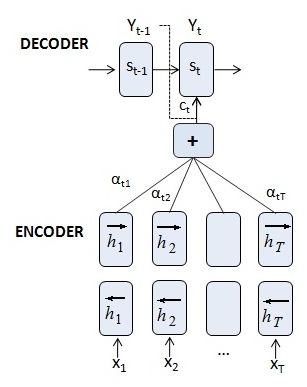} 
\caption{An encoder-decoder network with attention.}
\label{encatt}
\end{figure}

\noindent In this attention mechanism, the query is the state $\boldsymbol{s}_{t-1}$ and the key and the value are the hidden states $\boldsymbol{h}_{i}$. The score measures how well the input at position $i$ and the output at position $t$ match. $\alpha_{ti}$ are the weights that implement the attention mechanism, defining how much of each input hidden state should be considered when deciding the next state $\boldsymbol{s}_{t}$ and generating the output $\boldsymbol{y}_{t}$ (see Figure \ref{encatt}).

As we have described previously, the score function can be parametrized using different alignment models such as feedforward networks and the cosine similarity.

\begin{figure}[h]
\centering
\includegraphics[scale=0.6]{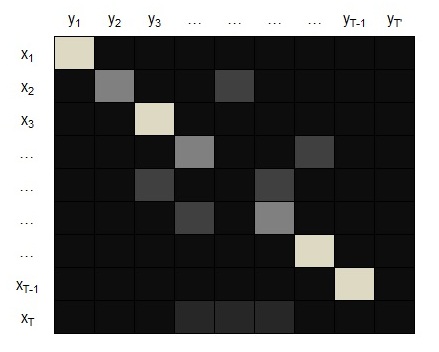} 
\caption{A matrix of alignment scores.}
\label{matrix}
\end{figure}

An example of a matrix of alignment scores can be seen in Figure \ref{matrix}. This matrix provides interpretability to the model since it allows to know which part (time-step) of the input is more important to the output.
 
\subsection{Other attention mechanisms and differentiable neural computers}

\label{sec32}

A variant of the attention mechanism is self-attention, in which the
attention component relates different positions of a single sequence in
order to compute a representation of the sequence. In this way, the keys, values and queries come from the same source. The mechanism can connect distant elements of the sequence more directly than using RNNs \cite{Tang2018SelfAtt}.

Another variant of attention are end-to-end memory networks \cite
{Suk2015EndToEndMN}, which we describe in Section \ref{sec422} and are neural networks with a recurrent attention model over an external memory. The model, trained end-to-end, outputs an answer based on a set of inputs $x_{1},x_{2},...,x_{n}$ stored in a memory and a query.

Traditional computers are based on the von Neumann architecture which has
two basic components: the CPU (Central Processing Unit), which carries out
the program instructions, and the memory, which is accessed by the CPU to
perform write/read operations. In contrast, neural networks follow a hybrid
model where synaptic weights perform both processing and memory tasks.

Neural networks and deep learning models are good at mapping inputs to
outputs but are limited in their ability to use facts from previous events
and store useful information. Differentiable Neural Computers (DNCs) \cite%
{Graves2016DNC} try to overcome these shortcomings by combining neural networks
with an external read-write memory.

As described in \cite{Graves2016DNC}, a DNC is a neural network, called the
controller (playing the role of a differentiable CPU), with an external
memory, an $N\times W$ matrix. The DNC uses differentiable attention mechanisms to
define distributions (weightings) over the N rows and learn the importance
each row has in a read or write operation.

To select the most appropriate memory components during read/write
operations, a weighted sum $w(i)$ is used over the memory locations $%
i=1,...,N.$ The attention mechanism is used in three different ways: 
\begin{enumerate}[(i)]

\item Access content (read or write) based on similarity. 

\item Time ordered access (temporal links) to recover the sequences in the order in which they were written.

\item Dynamic memory allocation, where the DNC assigns and releases memory based on usage percentage.

\end{enumerate} 

At each time step, the DNC gets an input vector and emits an output vector that is a function of the combination of the input vector and the memories selected.

DNCs, by combining the following characteristics, have very promising
applications in complex tasks that require both perception and reasoning:

\begin{enumerate}[(i)]

\item The classical perception capability of neural networks.

\item Read and write capabilities based on content similarity and learned by
the model.

\item The use of previous knowledge to plan and reason.

\item End-to-end differentiability of the model.

\item Implementation using software packages with automatic differentiation
libraries such as PyTorch, Tensorflow or similar.
\end{enumerate}

\subsection{Meta-plasticity and continuous learning}

\label{sec33}

The combination of geometric modules (classical neural networks) with
algorithmic ones adds new learning capabilities to deep learning models. In
the previous sections we have seen that one way to improve the learning process
is by focusing on certain elements of the input or a memory and making this attention differentiable.

Another natural way to improve the process of learning is to incorporate
differentiable primitives that add flexibility and adaptability. A source of
inspiration is neuromodulators, which furnish the traditional synaptic transmission with new computational and processing capabilities \cite{Hernandez2018MultilayerAN}.

Unlike the continuous learning capabilities of animal brains, which
allow animals to adapt quickly to the experience, in neural networks, once
the training is completed, the parameters are fixed and the network stops
learning. To solve this issue, in \cite{Miconi2018DifferentiablePT} a
differentiable plasticity component is attached to the network that helps
previously-trained networks adapt to ongoing experience.

The process to introduce the differentiable plastic component in the network is as follows. The activation $y_j$ of neuron $j$ has a conventional fixed weight $w_{ij}$ and a plastic component $\alpha_{ij}H_{ij}(t)$, where $\alpha_{ij}$ is a structural parameter tuned during the training period and $H_{ij}(t)$ a
plastic component automatically updated as a function of ongoing inputs and
outputs. The equations for the activation of $y_j$ with learning rate $\eta$%
, as described in \cite{Miconi2018DifferentiablePT}, are:

\begin{equation}  \label{eq:diffplast1}
y_j=\tanh{\left \{ \sum_{i\in inputs}(w_{ij}+\alpha_{ij}H_{ij}(t))y_i\right
\} },
\end{equation}

\begin{equation}  \label{eq:diffplast2}
H_{ij}(t+1)=\eta y_i y_j+(1-\eta)H_{ij}(t).
\end{equation}

Then, during the initial training period, $w_{ij}$ and $\alpha_{ij}$ are
trained using gradient descent and after this period, the model keeps
learning from ongoing experience.

\section{Dynamical systems and differentiable programming}
\label{sec4}

\subsection{Modeling dynamical systems with neural networks}

\label{sec41}

Dynamical systems deal with time-evolutionary processes and their corresponding
systems of equations. At any given time, a dynamical system has a state that
can be represented by a point in a state space (manifold). The evolutionary
process of the dynamical system describes what future states follow from the
current state. This process can be deterministic, if its entire future is
uniquely determined by its current state, or non-deterministic otherwise \cite{Layek2015BookDS} (e.g., a random dynamical system \cite{Arnold2003RandomDS}). Furthermore, it can be a continuous-time process, represented by differential equations or, as in this paper, a discrete-time process, represented by difference equations or maps. Thus,

\begin{equation}  \label{eq:dynamicalsystem1}
\boldsymbol{h}_t = f(\boldsymbol{h}_{t-1};\boldsymbol{\theta})
\end{equation}
\noindent for autonomous discrete-time deterministic dynamical systems with parameters $\boldsymbol{\theta}$, and

\begin{equation}  \label{eq:dynamicalsystem2}
\boldsymbol{h}_t = f(\boldsymbol{h}_{t-1},\boldsymbol{x}_t;\boldsymbol{\theta})
\end{equation}

\noindent for non-autonomous discrete-time deterministic dynamical systems driven by an external input $\boldsymbol{x}_t$.

Dynamical systems have important applications in physics, chemistry,
economics, engineering, biology and medicine \cite{Jackson2015BookDSBio}.
They are relevant even in day-to-day phenomena with great social impact such as
tsunami warning, earth temperature analysis and financial markets
prediction.

Dynamical systems that contain a very large number of variables
interacting with each other in non-trivial ways are sometimes called complex (dynamical) systems \cite%
{Gross2008CAS}. Their behaviour is intrinsically difficult to model due to
the dependencies and interactions between their parts and they have emergence
properties arising from these interactions such as adaptation, evolution,
learning, etc.

Here we consider discrete-time, deterministic and non-autonomous (i.e., the time evolution depending also on exogenous variables) dynamical systems as well as the more general complex systems. Specifically, the dynamical systems of interest range from systems of difference equations with multiple time delays to systems with a dynamic (i.e., time-changing) interdependence between time steps. Notice that the former ones may be rewritten as higher dimensional systems with time delay 1. 

On the other hand, in recent years deep learning models have been very
successful in performing various tasks such as image recognition, machine
translation, game playing, etc. When the amount of training data is sufficient 
and the distribution that generates the real data is the same as the distribution of the training data, these models perform extremely well and approximate the input-output relation.

In view of the importance of dynamical systems for modeling physical, biological
and social phenomena, there is a growing interest in applying deep
learning techniques to dynamical systems. This can be done in different contexts, such as:

\begin{enumerate}[(i)]

\item Modeling dynamical systems with known structure and equations but
non-analytical or complex solutions \cite{Pan2018LongTimePM}.

\item Modeling dynamical systems without knowledge of the underlying governing equations \cite%
{Düben2018GlobalWeather, Chakraborty1992ForecastingTB}. In this regard, let us mention  that commercial initiatives are emerging that combine large
amounts of meteorological data with deep learning models to improve weather
predictions.

\item Modeling dynamical systems with partial or noisy data \cite%
{Yeo2019DLNoisyDS}.
\end{enumerate}

A key aspect in modelling dynamical systems is temporal dependence. There
are two ways to introduce it into a neural network \cite%
{Narendra1990DSystemsNN}:

\begin{enumerate}[(i)]

\item A classical feedforward neural network with time delayed states in the
inputs but perhaps with an unnecessary increase in the number of parameters.

\item A recurrent neural network (RNN) which, as shown in Equations (\ref{eq:rnn1}) and (\ref{eq:rnn2}), has a temporal recurrence that makes it appropriate for modelling discrete dynamical systems of the form given in Equations (\ref{eq:dynamicalsystem1}) and (\ref{eq:dynamicalsystem2}). 
\end{enumerate}

Thus, RNNs, specially designed for sequence modelling \cite%
{chang2018antisymmetricrnn}, seem the ideal candidates to model, analyze and predict dynamical systems in the broad sense used in this tutorial. The temporal recurrence of RNNs, theoretically, allows to model and identify dynamical systems described with equations with any temporal dependence.

To learn chaotic dynamics, recurrent radial basis function (RBF) networks \cite%
{Miyosi1995ChaoticRBF} and evolutionary algorithms that generate RNNs have been proposed \cite{SatoEvolutionaryRNN}. "Nonlinear
Autoregressive model with exogenous input" (NARX) \cite%
{DiaconescuNARXChaoticSeries} and boosted RNNs \cite{AssaadChaoticSeriesRNN}
have been applied to predict chaotic time series.

\begin{figure}[h]
\centering
\includegraphics[scale=0.6]{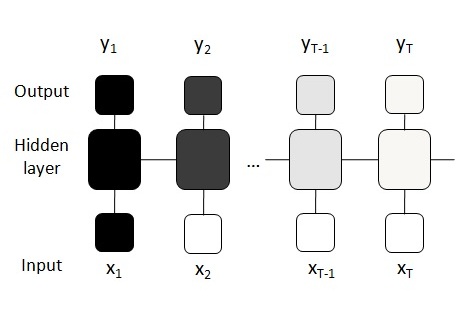} 
\caption{Vanishing gradient problem in RNNs. Information sensitivity decays over time forgetting the first input.}
\label{VanGrad}
\end{figure}

However, a difficulty with RNNs is the vanishing gradient problem \cite%
{Bengio1994GradientDifficult}. RNNs are trained by unfolding them into deep
feedforward networks, creating a new layer for each time step of the input
sequence. When backpropagation computes the gradient by the chain rule, this
gradient vanishes as the number of time-steps increases. As a result, for long
input-output sequences, as depicted in Figure \ref{VanGrad}, RNNs have trouble modelling long-term dependencies, that is, relationships between elements that are separated by large periods of time.

To overcome this problem, LSTMs were proposed. LSTMs have an advantage over basic RNNs due to their relative insensitivity to temporal delays and, therefore, are appropriate for modeling and making predictions based on time series whenever there exist temporary dependencies of unknown duration. With the appropriate number of hidden units and activation functions \cite{Yu2017ConceptLSTM}, LSTMs can model and identify any non-linear dynamical system of the form: 

\begin{equation}  \label{eq:dynamical1}
\boldsymbol{h}_{t}=f(\boldsymbol{x}_{t},...,\boldsymbol{x}_{t-T},\boldsymbol{h}_{t-1},...,\boldsymbol{h}_{t-T}),
\end{equation}
\begin{equation}  \label{eq:dynamical2}
\boldsymbol{y}_{t}=g(\boldsymbol{h}_{t}),
\end{equation}

\noindent $f$ and $g$ are the state and output functions while $\boldsymbol{x}_{t}$, $\boldsymbol{h}_{t}$ and $\boldsymbol{y}_{t}$ are the system input, state and output.

LSTMs have succeeded in various applications to dynamical systems
such as model identification and time series prediction \cite%
{WangModelIdLSTM, Yu2017ConceptLSTM, Li2018LSTMTourismFlow}.

An also remarkable application of the LSTM has been machine translation \cite{Sutskever2014SequenceTS, Cho2014EncDec}, using the encoder-decoder architecture described in Section \ref{sec311}.

However, as we have seen, the decoder possibly does not take into account the first elements of the input sequence because the encoder compresses all the information of the input sequence in a fixed-length vector. Then, the performance of encoder-decoder networks degrades rapidly as the length of
input sequence increases and this can be a problem in time series analysis,
where predictions are based upon a long segment of the series.

\begin{figure}[h]
\centering
\includegraphics[scale=0.48]{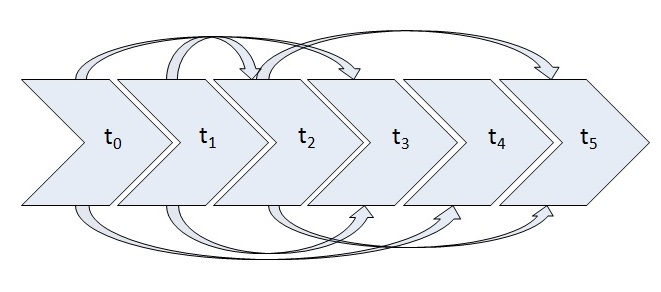} 
\caption{Temporal interdependencies in a dynamical system.}
\label{depend}
\end{figure}

Furthermore, as depicted in Figure \ref{depend}, a complex dynamic may feature interdependencies between time steps that vary with time. In this situation, the equation that defines the temporal evolution may change at each $t \in 1,...,T$. For these dynamical systems, adding an attention module like the one described in Equation \ref{eq:att1} can help model such time-changing interdependencies.

\subsection{Improving dynamical systems with differentiable programming}

\label{sec42}

Deep learning models together with graphic processors and large
amounts of data have improved the modeling of dynamical systems but this has some limitations such as those mentioned in the previous section. The combination of neural networks with new differentiable algorithmic modules is expected to overcome some of those shortcomings and offer new opportunities and applications.

In the next three subsections we illustrate with examples the kind of applications of differentiable programming to dynamical systems we have in mind, namely: implementations of attention mechanisms, memory networks, scientific simulations and modeling in physics.

\subsubsection{Attention mechanisms in dynamical systems}

\label{sec421}

In the previous sections we have described the attention mechanism, which allows a task to be guided by a set of elements of the input or memory source. When applying this mechanism to dynamical systems modeling or prediction, it is necessary to decide the following aspects:

\begin{enumerate}[(i)]
\item In which phase or phases of the model should the attention mechanism be introduced?

\item What dimension is the mechanism going to focus on? Temporal, spatial, etc.

\item What parts of the system will correspond to the query, the key and the value?

\end{enumerate}

One option, which is also quite illustrative, is to use a dual-stage attention, an encoder with input attention and a decoder with temporal attention, as pointed out in \cite{Quin2017DualAttention}. 

Here we describe this option, in which the first stage extracts the relevant input features and the second selects the relevant time steps of the model. In many dynamical systems there are long term dependencies between time steps and these dependencies can be dynamic, i.e., time-changing. In these cases, attention mechanisms learn to focus on the most relevant parts of the system input or state.

$X=(\boldsymbol{x}_{1},\boldsymbol{x}_{2},...,\boldsymbol{x}_{T})$ with $\boldsymbol{x}_{t}\in  \mathbb{R}^n$ represents the input sequence. $T$ is the length of the time interval and $n$ the number of input features or dimensions. At each time step $t$, $\boldsymbol{x}_{t}=({x}_{t}^{1},{x}_{t}^{2},...,{x}_{t}^{n})$. 

\vspace{5 mm}
\noindent \textbf{Encoder with input attention}

The encoder, given an input sequence $X$, maps $\boldsymbol{u}_{t}$ to

\begin{equation}  \label{eq:DualEnc1}
\boldsymbol{h}_{t}=f_{1}(\boldsymbol{h}_{t-1},\boldsymbol{u}_{t}),
\end{equation}

\noindent where $\boldsymbol{h}_{t} \in \mathbb{R}^m$ is the hidden state of the
encoder at time $t$, $m$ is the size of the hidden state and $f_{1}$ is an
RNN (or any of its variants). $\boldsymbol{x}_{t}$ is replaced by  $\boldsymbol{u}_{t}$, which adaptively selects the relevant input features with 

\begin{equation}  \label{eq:DualEnc2}
\boldsymbol{u}_{t}=(\alpha_{t}^{1}{x}_{t}^1,\alpha_{t}^{2}{x}_{t}^2,...,\alpha_{t}^{n}{x}_{t}^{n}).
\end{equation}

$\alpha_{t}^{k}$ is the attention weight measuring the importance of the $k$ input feature at time $t$ and is computed by

\begin{equation}  \label{eq:DualEnc3}
\alpha_{t}^{k}=\frac{\exp(score(\boldsymbol{h}_{t-1},\boldsymbol{x}^{k}))}{\sum_{i=1}^{T}\exp(score(\boldsymbol{h}_{t-1},\boldsymbol{x}^{i}))},
\end{equation}

\noindent where $\boldsymbol{x}^{k}=({x}_{1}^{k},{x}_{2}^{k},...,{x}_{T}^{k})$ is the $k$ input feature series and the score function can be computed using a feedforward neural network, a cosine similarity measure or other similarity functions. 

\begin{figure}[h]
\centering
\includegraphics[scale=0.53]{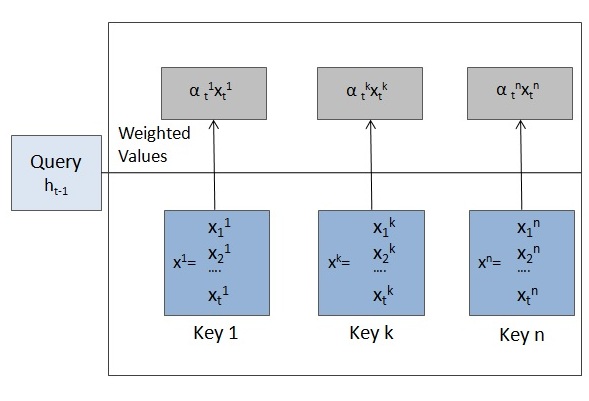} 
\caption{Diagram of the input attention mechanism.}
\label{InputStage}
\end{figure}

Then, this first attention stage extracts the relevant input features, as seen in Figure \ref{InputStage} with the corresponding query, keys and values.

\vspace{5 mm}
\noindent \textbf{Decoder with temporal attention}

Similar to the attention decoder described in Section \ref{sec311}, the decoder has an output

\begin{equation}  \label{eq:DualDecod1}
\boldsymbol{y}_{t}=f_{2}(\boldsymbol{s}_{t-1},\boldsymbol{y}_{t-1},\boldsymbol{c}_{t})
\end{equation}

\noindent for $t=1,...,T^{\prime }$. $f_{2}$ is an RNN (or any of its variants) with an additional linear or softmax layer, and the input is a concatenation of $\boldsymbol{y}_{t-1}$ with the context vector $\boldsymbol{c}_{t}$, which is a sum of hidden states of the input sequence weighted by alignment scores:

\begin{equation}  \label{eq:DualDecod2}
\boldsymbol{c}_{t}=\sum_{i=1}^{T}\beta_{t}^{i}\boldsymbol{h}_{i}.
\end{equation}

The weight $\beta_{t}^{i}$ of each state $\boldsymbol{h}_{i}$ is computed using the similarity function, $score(\boldsymbol{s}_{t-1},\boldsymbol{h}_{i})$, and applying a softmax function, as described in Section \ref{sec311}.

\begin{figure}[h]
\centering
\includegraphics[scale=0.53]{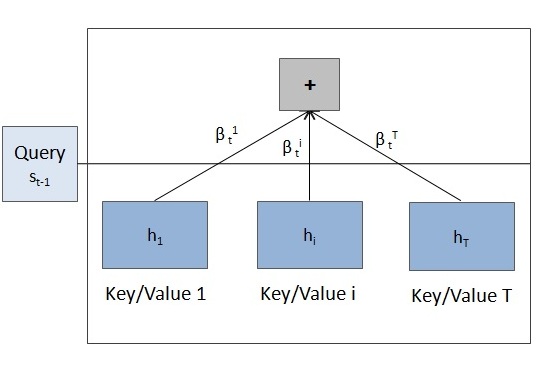} 
\caption{Diagram of the input attention mechanism.}
\label{TempStage}
\end{figure}

This second attention stage selects the relevant time steps, as shown in Figure \ref{TempStage} with the corresponding query, keys and values.

\vspace{5 mm}
\noindent \textbf{Further remarks}

In \cite{Quin2017DualAttention}, the authors define this dual-stage attention RNN and show that the model outperforms a classical model in time series prediction.

In \cite{Hollis2018LSTMvsAttention}, a comparison is made between LSTMs and
attention mechanisms for financial time series forecasting. It is shown there that an LSTM with attention perform better than stand-alone LSTMs.

A temporal attention layer is used in \cite{Vinayavekhin2018AttentonTimeSeries} to select relevant information and to provide model interpretability, an essential feature to understand deep learning models. Interpretability is further studied in detail in \cite{Serrano2019AttInterpret}, concluding that attention weights partially reflect the impact of the input elements on model prediction.

Despite the theoretical advantages and some achievements, further studies are needed to verify the benefits of the attention mechanism over traditional networks.

\subsubsection{Memory networks}

\label{sec422}

Memory networks allow long-term dependencies in sequential data to be learned thanks to an external memory component. Instead of taking into account only the most recent states, memory networks consider the entire list of entries or states.

Here we define one possible application of memory networks to dynamical systems, following an approach based on \cite{Suk2015EndToEndMN}. We are given a time series of historical data $\boldsymbol{n}_{1},...,\boldsymbol{n}_{T^\prime}$ with $\boldsymbol{n}_{i}\in  \mathbb{R}^n$ and the input series $\boldsymbol{x}_{1},...,\boldsymbol{x}_{T}$ with $\boldsymbol{x}_t \in \mathbb{R}^n$ the current input, which is the query.

The set $\{\boldsymbol{n}_{i}\}$ are converted into memory vectors $\{\boldsymbol{m}_{i}\}$ and output vectors $\{\boldsymbol{c}_{i}\}$ of dimension $d$. The query $\boldsymbol{x}_t$ is also transformed to obtain a internal state $\boldsymbol{u}_t$ of dimension $d$. These transformations correspond to a linear transformation: $A\boldsymbol{n}_{i}=\boldsymbol{m}_{i},B\boldsymbol{n}_{i}=\boldsymbol{c}_{i},C\boldsymbol{x}_t=\boldsymbol{u}_t$, being $A,B,C$ parameterizable matrices.

A match between $\boldsymbol{u}_t$ and each memory vector $\boldsymbol{m}_{i}$ is computed by taking the inner product followed by a softmax function:

\begin{equation}  \label{eq:MemNetwork1}
p_t^i=Softmax(\boldsymbol{u}_t^{T}\boldsymbol{m}_{i}).
\end{equation} 

The final vector from the memory, $\boldsymbol{o}_t$, is a weighted sum over the transformed inputs $\{\boldsymbol{c}_{i}\}$:

\begin{equation}  \label{eq:MemNetwork2}
\boldsymbol{o}_t=\sum_{j}p_t^i\boldsymbol{c}_{i}.
\end{equation}

To generate the final prediction $\boldsymbol{y}_t$, a linear layer is applied to the sum of the output vector $\boldsymbol{o}_t$ and the transformed input $\boldsymbol{u}_t$ and to the previous output $\boldsymbol{y}_{t-1}$:

\begin{equation}  \label{eq:MemNetwork3}
\boldsymbol{y}_t=W^1(\boldsymbol{o}_t+\boldsymbol{u}_t)+W^2\boldsymbol{y}_{t-1}
\end{equation}

This model is differentiable end-to-end by learning the matrices (the final matrices $W^i$ ant the three transformation matrices $A,B$ and $C$) to minimize the prediction error.

In \cite{Chang2018AMB} the authors propose a similar model based on memory networks with a memory component, three encoders and an autoregressive component for multivariate time-series forecasting. Compared to non-memory RNN models, their model is better at modeling and capturing long-term dependencies and, moreover, it is interpretable.

Taking advantage of the highlighted capabilities of Differentiable Neural
Computers (DNCs), an enhanced DNC for electroencephalogram (EEG)
data analysis is proposed in \cite{Ming2018EEGDataDNC}. By replacing the LSTM network controller with a recurrent convolutional network, the
potential of DNCs in EEG signal processing is convincingly demonstrated.

\subsubsection{Scientific simulation and physical modeling}

\label{sec423}

Scientific modeling, as pointed out in \cite{Rackauckas2019DiffEq}, has traditionally employed three approaches: 
\begin{enumerate}[(i)]

\item Direct modeling, if the exact function that relates input and output is known.

\item Using a machine learning model. As we have mentioned, neural networks are universal approximators.

\item Using a differential equation if some structure of the problem is known. For example, if the rate of change of the unknown function is a function of the physical variables.    

\end{enumerate}

Machine learning models have to learn the input-output transformation from scratch and need a lot of data. One way to make them more efficient is to combine them with a differentiable component suited to a specific problem. This component allows specific prior knowledge to be incorporated into deep learning models and can be a differentiable physical model or a differentiable ODE (ordinary differential equation) solver.

\begin{enumerate}[(i)]

\item \textit{Differentiable physical models}.

Differentiable plasticity, as described in Section \ref{sec33}, can be applied to deep learning models of dynamical systems in order to help them adapt to ongoing data and experience. 

As done in \cite{Miconi2018DifferentiablePT}, the plasticity component described in Equations \ref{eq:diffplast1} and \ref{eq:diffplast2}, can be introduced in some layers of the deep learning architecture. In this way, the model can continuously learn because the plastic component is updated by neural activity.

DiffTaichi, a differentiable programming language for building differentiable
physical simulations, is proposed in \cite{Hu2019DiffTaichi}, integrating a
neural network controller with a physical simulation module.

A differentiable physics engine is presented in \cite{Belbute2018DiffPhysics}. The system simulates rigid body dynamics and can be integrated in an
end-to-end differentiable deep learning model for learning the physical
parameters.

\item \textit{Differentiable ODE solvers}.

As described in \cite{Rackauckas2019DiffEq}, an ODE can be embedded into a deep learning model. For example, the Euler method takes in the derivative function and the initial values and outputs the approximated solution. The derivative function could be a neural network. 

This solver is differentiable and can be integrated into a lager model that can be optimized using gradient descent.

In \cite{Innes2019Zygote} a differentiable model of a trebuchet is described. In a classical trebuchet model, the parameters (the mass of the counterweight and the angle of release) are fed into an ODE solver that calculates the distance, which is compared with the target distance. 

In the extended model, a neural network is introduced. The network takes two inputs, the target distance and the current wind speed, and outputs the trebuchet parameters, which are fed into the simulator to calculate the distance. This distance is compared with the target distance and the error is back-propagated through the entire model to optimize the parameters of the network. Then, the neural network is optimized so that the model can achieve any target distance. Using this extended model is faster than optimizing only the trebuchet.

This type of applications shows how combining differentiable ODE solvers and deep learning models allows to incorporate previous structure to the problem and makes the learning process more efficient.    

We may conclude that combining scientific computing and differentiable components will open new avenues in the coming years.

\end{enumerate}

\section{Conclusions and future directions}

\label{sec6}

Differentiable programming is the use of new differentiable components beyond classical neural networks. This generalization of deep learning allows to have data parametrizable architectures instead of pre-fixed ones and new learning capabilities such as reasoning, attention and memory.

The first models created under this new paradigm, such as attention mechanisms, differentiable neural computers and memory networks, are already having a great impact on natural language processing.

These new models and differentiable programming are also beginning to improve machine learning applications to dynamical systems. As we have seen, these models improve the capabilities of RNNs and LSTMs in identification, modeling and prediction of dynamical systems. They even add a necessary feature in machine learning such as interpretability.

However, this is an emerging field and further research is needed in several directions. To mention a few:

\begin{enumerate}[(i)]

\item More comparative studies between attention mechanisms and LSTMs in predicting dynamical systems. 

\item Use of self-attention and its possible applications to dynamical systems. 

\item As with RNNs, a theoretical analysis (e.g., in the framework of dynamical systems) of attention and memory networks.

\item Clear guidelines so that scientists without advanced knowledge of machine learning can use new differentiable models in computational simulations.

\end{enumerate}

\medskip

\textit{Acknowledgments.} This work was financially supported by the Spanish Ministry of
Science, Innovation and Universities, grant MTM2016-74921-P (AEI/FEDER, EU).

\appendix




\bibliographystyle{elsarticle-num}
\bibliography{sample}

\end{document}